# Parametric Identification of a Functional-Structural Tree Growth Model and Application to Beech Trees (*Fagus sylvatica*)

*Véronique Letort*[A,*]*, Paul-Henry Cournède*[A]*, Amélie Mathieu*[A]*, Philippe de Reffye*[B]*, Thiéry Constant*[C]

[A]Ecole Centrale of Paris, Laboratoire de Mathématiques Appliquées aux Systèmes, F-92295 Châtenay-Malabry cedex, France
[B]Cirad-Amis, UMR AMAP, TA 40/01 Ave Agropolis, F-34398 Montpellier cedex 5, France and INRIA-Saclay, Parc Orsay Université, F-91893 Orsay cedex, France
[C]LERFOB UMR INRA-ENGREF No. 1092, Wood Quality Research Team, INRA Research Centre of Nancy, F-54280 Champenoux, France.

[*]Corresponding author. Email: letort@mas.ecp.fr

**Abstract.** Functional-structural models provide detailed representations of tree growth and their application to forestry seems full of prospects. However, due to the complexity of tree architecture, parametric identification of such models remains a critical issue. We present the GreenLab approach for modeling tree growth. It simulates tree plasticity in response to changes of their internal level of trophic competition, especially regarding topological development and cambial growth. The model includes a simplified representation of tree architecture, based on a species-specific description of branching patterns. We study whether those simplifications allow enough flexibility to reproduce with the same set of parameters the growth of two observed understorey Beech trees (*Fagus sylvatica*, L.) of different ages and in different environmental conditions. The parametric identification of the model is global, *i.e.* all parameters are estimated simultaneously, potentially providing a better description of interactions between sub-processes. As a result, the source-sink dynamics throughout tree development is retrieved. Simulated and measured trees were compared for their trunk profiles (fresh masses and dimensions of every growth units, ring diameters at different heights) and for the compartment masses of their order 2 branches. Possible improvements of this method by including topological criteria are discussed.

## Introduction

In tree growth simulation, plant architecture can be described at different levels of simplification depending on the objectives. In forestry models, the usual way to simplify crown shape is to consider only its exterior surface and to fit the equation of a species-dependent geometrical surface (Pretzsch *et al* 2002). But some processes (e.g. light interception, biomass allocation) and some descriptive variables (e.g. crown surface, ring, profiles, competition index) could be assessed more mechanistically if a more detailed description of tree architecture was taken into account, especially in heterogeneous stands where species-specific architectural patterns and relative positions of neighbour trees play a predominant role in determining the effects of competition. In that context, a potentially interesting framework is provided by the recent emergence of functional-structural models (FSMs) whose objective is to simulate interactively the architectural development of trees and their physiological functioning (Sievänen *et al.* 2000; Prusinkiewicz 2004).

Some of the existing FSMs provide a fine description of the source-sink interactions at organ level. The biomass production of each leaf is usually computed on a hourly or daily basis according to the amount of intercepted radiations. Dynamics of allocation patterns are defined with various degrees of details; for instance using the transport-resistance model in L-Peach (Grossman and DeJong 1994; Allen *et al*. 2005), with a matrix of measured allocation coefficients for each source-sink pair in ECOPHYS (Rauscher *et al.* 1990), using fitted regression equations for allocation fractions in Zhang *et al.* (1994) or including a function representing the source-sink distance in SIMWAL (Balandier *et al.* 2000). These models are mainly dedicated to simulating the growth of fruit trees or young trees for periods ranging from several months to a few years.

Other models aim at simulating the growth of forest trees with periods of interest that can largely exceed one hundred years. They generally incorporate a more limited set of processes and variables with a time step of one year. For instance, the LIGNUM model (Perttunen *et al.* 1996; Perttunen and Sievänen 2005) or the model of Sterck *et al.* (2005) account for the local effects of light interception on biomass production, allocation and architectural development. Rather than a detailed reproduction of the tree structure, their models allow the exploration of the global effects of model assumptions and parameter variations (Sterck and Schieving 2007). These models have to face the problem of simulation performance when dealing with high numbers of organs (artificially limited to 5000 in Sterck and Schieving 2007) or when considering the growth of heterogeneous stands within a reasonable time.

However, despite the research effort involved in their development, the use of these models remains confined to research and teaching contexts; forestry applications are scarce (Le Roux *et al.* 2001). Besides the problem of simulation efficiency, a reason could be that parametric identification (that is to say the estimation of model parameters, *cf.* Walter and Pronzato, 1994) and validation of such models are still critical points. Indeed, although these models provide very fine descriptions of trees at organ scale, experimental data are not usually available at the same scale due to high topological complexity and tremendous number of organs in trees: only global, aggregated or sampled measurements are reasonable to expect. For instance, LIGNUM has been parameterized for Scots pine (Perttunen *et al.* 1996) and Sugar maple (Perttunen *et al.* 2001) considering independently the physiological processes involved and comparing the model outputs to more aggregated data (e.g. taper curve, height, leaf area).

In this context, the objective of this paper is to present the GreenLab approach for functional-structural modeling of tree growth and the methodology for its parametric identification. GreenLab simulates biomass production and allocation using a source-sink model. The time step (called growth cycle) is equal to one year for temperate trees without polycyclism. Hereafter the term biomass refers to fresh matter (we refer to Louarn *et al.* 2008 regarding the relationships between fresh and dry matter content in the context of GreenLab source-sink dynamics). It provides a topological description of plant architecture but with simplified rules, which could help bridging the gap between functional-structural modeling and forestry applications. These simplified rules rely on the principle of structure factorization as defined in Yan *et al.* (2004) and Cournède *et al.* (2006). It takes advantage of repetitions in tree architecture (Barthélémy and Caraglio 2007) and the description of tree structure can be written under the form of a single recursive equation. As it is consistent with the multi-scale decomposition of tree architecture proposed by Godin and Caraglio (1998), some species-specific architectural trends identified thanks to botanical analyses can be taken into account. It also allows the formalization of the model as a dynamic system for which classical methods of parametric identification can be applied. Contrary to other models whose parameters must be estimated on independent data sets for each functional sub-process, identification of GreenLab is global, that is to say that all parameters are estimated simultaneously, potentially providing a better description of interactions between sub-processes. It is thus possible to track back the dynamics of source-sink relationships between organs throughout tree growth.

Plasticity is an important feature in tree growth. It is modeled in GreenLab through tree adaptation to internal states of trophic competition, especially regarding topological development and cambial growth (Mathieu *et al.*, 2006; Mathieu 2006). Trophic competition is described in terms of source-sink dynamics as the ratio of available biomass to the demand of all sinks. For instance, lateral buds will develop into new branches only when the ratio of biomass to demand exceeds a threshold value (see also the Ecomeristem model in Dingkuhn *et al.* 2006; Luquet *et al.* 2006; Luquet *et al.* 2007). Tree topological development adapts to its available resources and the level of competition between organs: the simulated tree reacts as a self-regulating system. The step of parametric identification presented in this paper explores the possibility of extracting the values of these thresholds from experimental data and of using such a self-regulating system to generate coherent tree branching patterns.

In a first part, we present the specificities of GreenLab for trees and we describe the parametric identification procedure from experimental data. In a second part, we show an application of the method to two common beech trees (*Fagus Sylvatica*, L.) grown in a forest near Nancy (North-Eastern part of France). Although the two measured trees were very different regarding their ages and dimensions, the same set of model parameters is used to fit the data of both trees; only one control variable representing the relative level of environmental pressure on each tree differs.

## Materials and methods

For a detailed presentation of the concepts underlying the GreenLab functional-structural model, we refer to Yan *et al.* (2004) and Cournède *et al.* (2006). We only recall here its basic principles and mostly focus on its specificities for the modeling of tree growth.

### *A source-sink model*

Let us consider an individual tree. Its biomass production and allocation are driven by a source-sink model. Since few biomass data are available concerning tree root system, we assume a constant proportion of allocation to the underground growth at each cycle. Hereafter, *Q(n)* denotes the biomass production allocated to the aerial part at growth cycle *n*, once subtracted the amount of assimilates required for respiration processes. The net production is calculated according to light interception by the crown, expressed using classical Beer's law. Several adaptations of Beer's law were proposed for trees (Wang 2003; Sinoquet *et al.* 2005). Here we use the following simplified form:

$$Q(n) = PAR(n) \cdot RUE \cdot Sp(n) \cdot \left(1 - e^{-k \frac{S(n)}{Sp(n)}}\right) \quad (1)$$

where *PAR(n)* is the incident photosynthetically active radiation during growth cycle *n* and *RUE* is the radiation use efficiency. *k* is the extinction coefficient for radiations propagating vertically in the crown. *S(n)* is the total blade area of the plant at growth cycle *n*. *Sp(n)* is a parameter representing a characteristic surface area, related to the projection surface of the crown and modulated by the effects of competition with neighbours (Cournède *et al.* 2008). The long-term development of tree crown and the step by step metamorphosis of the tree architecture make the ratio of *Sp(n)* to *S(n)* variable. In GreenLab, a simple allometric relationship was chosen to compute the evolution of the projected surface according to the blade area:

$$Sp(n) = Sp_0 \cdot \left(\frac{S(n)}{Sp_0}\right)^{\alpha} \quad (2)$$

Parameters $Sp_0$ and $\alpha$ are unknown parameters that must be estimated by model inversion so that the variation of *Sp* induces the measured value of biomass production.

At each growth cycle, biomass is allocated to two distinct growth compartments: primary growth ($Q_s$) and secondary growth ($Q_r$) as represented in Fig. 1.

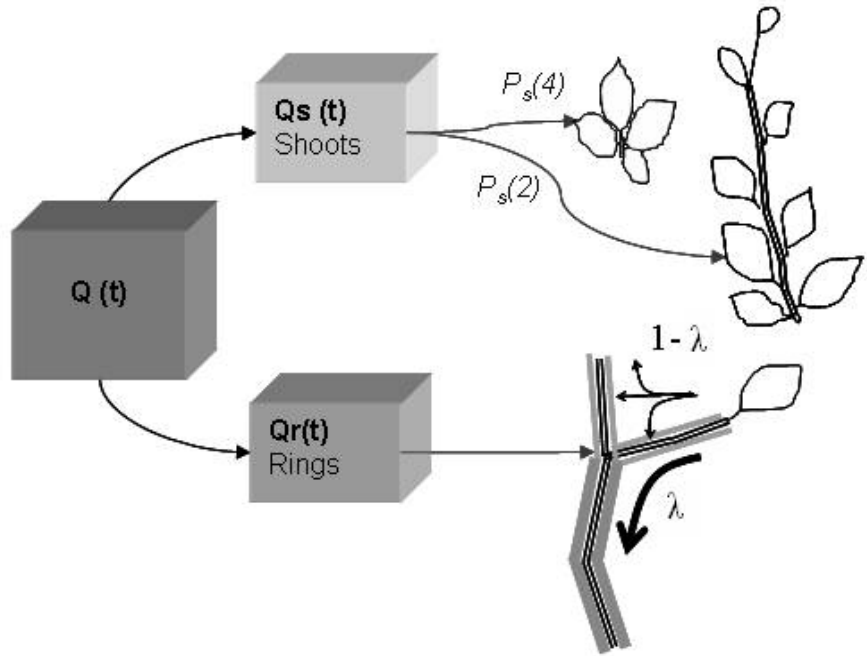

**Fig. 1**. *Model of biomass allocation. In a first step, biomass is allocated to compartments at the whole plant scale (new shoots and rings). In a second step, it is partitioned between all active metamers, according to their respective sink strengths ($P_s(2)$ and $P_s(4)$ for long and short shoots). The coefficient λ in [0;1] drives the relative influence on ring growth of*

*the amount of leaf area above the metamer in the topological structure.*
Primary growth corresponds to the emergence and the expansion of new shoots. It occurs at the beginning of the growth cycle and has a short duration (around two weeks) compared to that of the growth cycle. Secondary growth corresponds to the radial biomass increment of existing metamers and lasts during the whole growth cycle. Thus the growth of new shoots depends on prior-year conditions whereas the amount of cambial growth depends on the current year conditions (Duff and Nolan 1953, White 1993). Biomass is allocated to these two compartments proportionally to their respective demands, as described in the next two paragraphs.

### *New shoots and topological plasticity*

The demand of new shoots appearing at growth cycle *n+1* is calculated at the end of the previous growth cycle and is thus denoted $D_s(n)$. It is calculated as the sum of all potential shoot demands. In GreenLab, new shoots can be sorted into categories based on the botanical notion of physiological age (PA), as defined in Barthélémy *et al.* (1997). It characterizes the stage of a meristem in its potential ontogenetic evolution. Usually, no more than five physiological ages are necessary to characterize all the types of axes in a tree (Sabatier and Barthélémy 1999). A set of parameters values can be attached to each PA-based category. We have:

$$D_s(n) = \sum_{k=1}^{PA_m} N_{bud}(k,n) \cdot P_s(k) \quad (3)$$

where $N_{bud}(k,n)$ represents the number of buds of physiological age *k* at cycle *n* and $P_s(k)$ is the sink of a shoot of physiological age *k*. The biomass allocated to a shoot of physiological age *k* that will appear at growth cycle *n+1* is thus:

$$q_s(k,n) = P_s(k) \cdot \frac{Q_s(n)}{D_s(n)} \quad (4)$$

Inside each shoot, the biomass $q_s(k,n)$ is divided into internodes and leaves with a constant ratio. The organ dimensions are computed from their biomass content using simple allometric relationships as in Guo *et al.* (2008).

The positions and the numbers of potential shoots ($N_{bud}(k,n)$) depend on the tree topology at cycle *n*. Growth units are decomposed into zones that are characterized by the production of lateral buds of their metamers, as introduced by Guédon *et al* (2001). Here we define a zone $Z^{ik}$ on a growth unit of physiological age *i* as the set of metamers bearing buds of physiological age *k* (Fig. 2). We note $M^{ik}$ the number of metamers in a zone $Z^{ik}$. The sequence of potential zones along growth units of each physiological age provides the tree default topology. It results from species-specific rules defining potential growth (*e.g.* maximal physiological age of buds) and branching patterns (*e.g.* maximal branching order, general structure of branching hierarchy). The number of metamer per zone and their lateral productions are updated at each growth cycle as functions of the ratio of plant biomass production to its demand, *Q/D* (Mathieu *et al.* 2006). It is hypothesized that this ratio characterizes the level of trophic competition inside the plant. A low value of that variable means a high level of competition between organs for biomass supply.

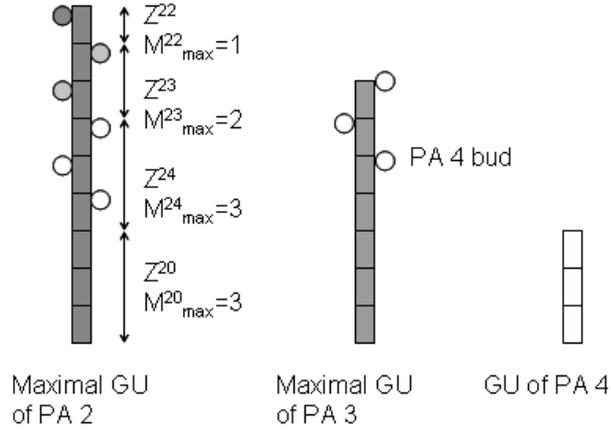

**Fig. 2**. *Default topology for growth units (GUs) of physiological ages 2, 3 and 4. The topology of GUs on the main stem (PA 1) is set from the target tree.*

At growth cycle *n*, the total number of new axes $A^{ik}(n)$ appearing on zones $Z^{ik}$ in the whole tree (*i.e.* the number of branches of PA *k* borne by metamers of PA *i*) depends on the number of existing positions potentially bearing that kind of axes $N^{ik}(n)$ and on the value of *Q/D*. It is given by the following equation:

$$A^{ik}(n) = \left[ N^{ik}(n) \cdot \left( A_1^{ik} + A_2^{ik} \cdot \frac{Q(n)}{D(n)} \right) \right] \quad (5)$$

where for a real number x, *[x]* represents its rounded value and $A^{ik}_1$ and $A^{ik}_2$ are model parameters. The $A^{ik}(n)$ axes are distributed on the $N^{ik}(n)$ positions under the following constraints: (*i*) growth units that appeared at the same growth cycle, with same physiological age and same rank along their mother axis bear the same number of axes, (*ii*) axes distribution begins from the oldest branches to the youngest ones (less ramified) and (*iii*) the distribution is as uniform as possible (all existing positions receive one axis before one position receive one more axis).

Similarly, the number of metamers $M^{ik}(n)$ per growth unit zone is determined according to the default topological rules and modulated by a function of *Q/D*:

$$M^{ik}(t) = \min\left( \left[ M_1^{ik} + M_2^{ik} \cdot \frac{Q(n)}{D(n)} \right], M_{\max}^{ik} \right) \quad (6)$$

where $M^{ik}_{max}$ is the maximal number of metamers in zone $Z^{ik}$, *i.e.* the maximal number of metamers of PA *i* that can potentially bear axes of PA *k*. $M^{ik}_1$ and $M^{ik}_2$ are model parameters. Note that the number of metamers per zone is calculated at growth unit level whereas the number of axes is determined at the whole-plant level.

The parameters $A^{ik}_1$, $A^{ik}_2$, $M^{ik}_1$ and $M^{ik}_2$ can be estimated from sampled observations or estimated by model inversion.

## *Cambial growth*

There are numerous evidences in the literature that cambial growth is more sensitive to changes in growth conditions than primary growth (Collet *et al.* 2001). Trees suffering from high competition can produce long primary shoots but usually narrow rings (Lanner 1985; Nicolini *et al.* 2001). Cournède *et al* (2008) showed that GreenLab mechanistically generates that behavior. Cambial growth is set to depend on the 'vigor' of the plant, which is associated to the main driving variable of the model, the ratio of plant biomass production to demand ($Q/D$) (Mathieu 2006). The ring demand is calculated as a power function of this ratio:

$$D_r(n) = P_r \cdot \left(\frac{Q(n)}{D(n)}\right)^\gamma \quad (7)$$

where $\gamma$ is an empirical parameter to be estimated. The proportion coefficient $P_r$ is the sink of the ring compartment. Thus the ring compartment plays a buffer role: the simulated tree invests more in secondary growth if the conditions are favorable.

Consequently, the global demand ($D_s + D_r$) of the tree is computed after solving the equation (8) for the unknown variable $D(n)$:

$$D(n) = D_s(n) + P_r \cdot \left(\frac{Q(n)}{D(n)}\right)^\gamma \quad (8)$$

This equation has a unique positive solution and can be numerically solved using the classical Newton method

In a second step (Fig. 1), the biomass globally allocated to the ring compartment is distributed to each metamer. In several models (reviewed in Le Roux *et al.* 2001), this step of biomass partitioning is done according to the pipe model principle introduced by Shinozaki *et al.* (1964): each leaf is associated to a woody pipe running to the stem base. This principle implies the Pressler rule that states that the annual ring surface is proportional to the amount of foliage above the ring position. However, Pouderoux *et al.* (2001) and Deleuze and Houllier (2002) pointed out several limitations of the Pressler rule. In order to investigate the influence of foliage biomass on ring growth, a flexible sub-model of allocation with two modes of partitioning is defined for biomass distribution to metamer rings in GreenLab: one is based on the Pressler rule and the other one is uniform allocation. In the first mode, the demand of a metamer for its cambial growth is calculated according to its secondary sink, its length and the surface of active leaves above it in the architecture. In the second mode, leaf surface is not considered and all leaves participate in equal proportions to the cambial growth of each metamer. These two modes can be mixed with a proportion coefficient $\lambda$ in [0;1], as represented in Fig. 1. More precisely, the set of equations driving the ring increment at growth cycle $n$ of a metamer of PA $k$, rank $p$ and chronological age $j$ can be written as follows:

$$\begin{cases} q_{rg}(k,p,j,n) = \left(\frac{1-\lambda}{D_{pool}(n)} + \frac{\lambda \cdot S_a(k,p,j,n)}{D_{pressler}(n)}\right) \cdot p_{rg}(k) \cdot l(k,j,n) \cdot Q_r(n) \\ D_{pressler}(n) = \sum_{k=1}^{PA_m} \sum_{p=1}^{n} \sum_{j=1_{k,i}}^{n} N_i(k,p,j,n) \cdot S_a(k,p,j,n) \cdot p_{rg}(k) \cdot l(k,j,n) \\ D_{pool}(n) = \sum_{k=1}^{PA_m} \sum_{i=1}^{n} N_i(k,j,n) \cdot p_{rg}(k) \cdot l(k,j,n) \end{cases} \quad (9)$$

where $S_a(k,p,j,n)$ denotes the surface of leaves above the considered metamer in the plant topology. It represents the number of leaves « seen » by the metamer. It is computed efficiently thanks to the plant decomposition into substructures (see Cournède *et al.* 2006). $N_i(k,p,j,n)$ is the total number of internodes of physiological age *k*, rank *p* and chronological age *j* in the plant at cycle *n*. $p_{rg}(k)$ is the linear sink of metamer of physiological age *k*. It is multiplied by the internode length *l(k,i,t)*. $D_{pressler}(n)$ and $D_{pool}(n)$ represent respectively the plant demand at cycle *n* for ring growth with the first and second modes. This model is more general than the Pressler rule that corresponds to the case $\lambda = 1$. The $\lambda$ coefficient is estimated from data and thus the fitting procedure can be used to assess the level of influence of the numbers and positions of leaves on the partitioning of ring biomass.

## *Measurement data on beech trees and calibration procedure*

Two common beech trees (*Fagus sylvatica* L.) growing as understorey trees were destructively measured in May 2006 from the natural stand of Champenoux, near Nancy (North-Eastern part of France). The two trees had not reached their sexual maturity at the measurement date: beech trees in dense stands generally reach sexual maturity after 60 to 80 years (Thiébaut and Puech 1984).

### *Beech topology*

The topological development of Beech tree follows the botanical rules of the Troll architectural model (Hallé and Oldeman 1970). The separations between annual shoots can be identified by morphological markers, namely the scars left by hibernal buds (Nicolini 1997). Annual shoots can be composed of several growth units separated by scars of ephemeral buds. In our observations, polycylic shoots were scarce, therefore this phenomenon was not considered, which is reasonable for understorey trees.

Beech development relies on two kinds of shoots: short shoots composed of a few metamers with negligible internodes and long shoots with elongated internodes (see Fig. 3).

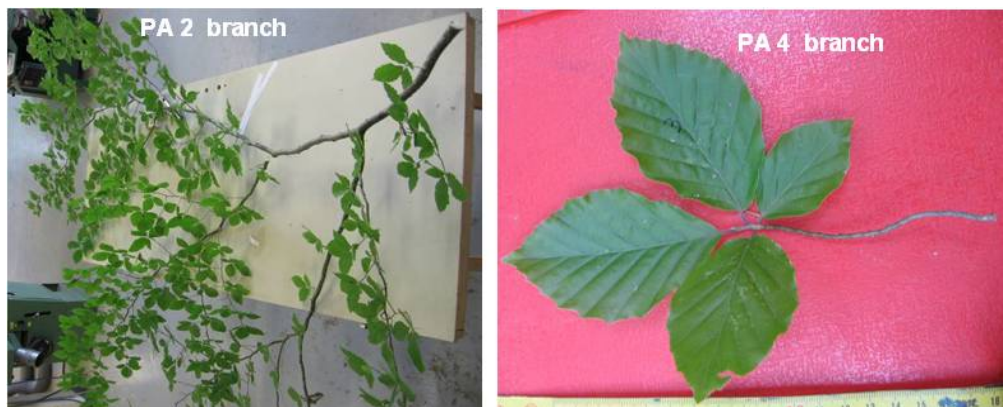

**Fig. 3**. *Examples of branches of physiological ages 2 and 4.*

Long shoots are dedicated to space exploration whereas short shoots are high contributors to plant production (Thiébaut and Puech 1984). We sorted the shoots into four PA-based categories. The trunk was attributed PA 1. Following Thiébaut and Puech (1984), we distinguished three classes of axes based on their elongation and ramification characteristics: short shoots bearing no branches (PA 4), long shoots bearing only short shoots (PA 3) and long ramified shoots (PA 2). The default topology for these three kinds of growth units is defined in Fig. 2. Growth units of PA 2 consist of four potential zones composed of metamers that can bear respectively no branches, axes of PA 4, PA 3 and even PA 2 (partial reiterations). Branches of PA 2 and PA 3 are both composed of long shoots with similar physiological characteristics but differ regarding their ramification parameters. The sequence of zones along growth units follows the botanical principle of acrotony, as described by Rauh for beech trees (Nicolini 1998).

*Direct measurements of model parameters*

The positions and numbers of metamers in the growth units of the main stem were recorded where the separation markers were still visible (which was not possible at the stem basis). The numbers and positions of branches were also observed, including scars indicating the positions of dead branches. To allow the comparison of the biomass compartment data between observations and simulations, the numbers and positions of branches on the trunk were imposed from the measurement. It means that the variations driven by the topological parameters (namely $A^{ik}_1$, $A^{ik}_2$, $M^{ik}_1$ and $M^{ik}_2$) were considered only for physiological ages higher than 1 ($i=2, 3$). Moreover, the values of some parameters of the model can be set thanks to the botanical knowledge on beech trees. For short shoots of PA 4, the number of metamers is relatively stable (Thiébaut and Puech 1984), usually 3 or 4. The constant value of 3 was uniformly chosen. It was also assumed that each zone $Z^{ik}$ has at least one metamer, except for the zone $Z^{22}$ (of PA 2 bearing lateral buds of PA 2). This is reasonable as growth units generally contain more than three metamers. By contrast, no minimal number of axes is set, as some growth unit zones may bear no lateral branches. These assumptions can be written as:

$$\forall\ i \in \{2,\ 3\},\ k \in \{0,\ 2,\ 3,\ 4\},\ \begin{cases} A^{ik}_1 = 1\ if\ (i,k) \neq (2,2),\ A^{ik}_1 = 0\ else \\ M^{ik}_1 = 0 \end{cases} \qquad (10)$$

Therefore only 10 topological parameters remained to be estimated (see section Results).

Other parameters, driving the physiological part of the model, can be set from direct observations. At random locations, blade surfaces were measured on samplings in order to have an estimation of the specific leaf weight (SLW). Nicolini and Chanson (1999) report that SLW also varies with the tree age because of ontogenic modifications of cell structures. Consequently, the same SLW evolution was set in the simulation for the two trees according to their age. Each sample consisted of blades belonging to the same branch (up to 50 blades per sample). Their fresh masses were measured then they were scanned and analyzed using the open source software ImageJ (ImageJ 1.36b, Wayne Rasband, National Institute of Health, USA).

These samples were also used to get the average ratios between the weight of one internode and of one leaf for the different kinds of shoots. These ratios provide the organ sinks in new shoots, *i.e.* the allometric coefficients for biomass partitioning between internode and leaf in each metamer.

The initial biomass ($Q_0$) was arbitrarily set to 1. As beech branches are plagiotropic (at least for young trees), most leaves are horizontally oriented so the coefficient $k$ for light interception is given the value 1 in Eqn 1. Mathematically speaking, the two factors *PAR(n)* and *RUE* in Eqn (1) can be gathered into one single factor, *V(n)*. As no measurement was done concerning this factor, it was assessed to be constant during the growth period but specific to the local environment of each tree. In our application, it was therefore considered as a parameter of the model (instead of a control factor) and its value was estimated from experimental biomass data. It is an indicator of the global level of environmental pressure on the tree growth.

Eventually Equation (1) is written as:

$$Q_t(n) = V_t \cdot Sp_0 \cdot \left(\frac{S_t(n)}{Sp_0}\right)^\alpha \cdot \left(1 - e^{-\left(\frac{S_t(n)}{Sp_0(n)}\right)^{1-\alpha}}\right) \quad (11)$$

Here the index $t$ ($t=1, 2$) is added to identify the variables that are specific to each tree, by contrast to the model parameters that are assumed to be species-specific and consequently take common values for both trees.

*Experimental data for parametric identification*

The experimental data are the target data to fit by the model, and the parameters are estimated by minimizing the distance (from a statistical point of view) between these target data and the corresponding simulations. In adequacy with the simplification level of the architectural description in our model, we chose to evaluate the adequacy between experimental observations and model simulations from the complete biomass profile of the trunk and biomass compartments on each branch of order two.

For each growth unit of the main stem, fresh mass, mean diameter and length were measured. At regular intervals, annual ring widths were recorded along four directions and averaged.

For every branch of order 2, its total fresh weigh, its main axis length and basal diameter were measured. All branches were considered (regardless of their basal diameters) but the method used to get their relative mass of wood and leaves were different depending on the complexity of the branch structure. For PA 4 and PA 3 branches, leaves were separated from internodes and weighted. Concerning PA 2 branches, the number of leaves was too high to get complete measurements. The following procedure was adopted to get an estimation of the wood and leaf weights: the basal part of the branch, bearing no leaves, was weighted separately and recorded as a woody part. For the rest of the branch, a representative branchlet was arbitrarily chosen. Its leaf and wood weights were measured. The same wood/leaf ratio was kept to get the compartment weights from the total weight of that part of the branch.

If several branches of same PA are located on the same GU, the average values for the leaf and internode compartments are put in the target. Of course, the cumulated values remain the same when summed on all branches.

To summarize, the following data were included in the target file: average internode length, diameter and fresh mass for each growth unit on the trunk, average ring diameters for some of those growth units, weights of leaf and wood compartments for each branch of order 2.

*Identification procedure*

These data were fitted with GreenLab using the DigiPlant software (Cournède *et al.* 2006). The objective criterion is calculated as the weighted sum of the square differences of target to simulated data. For some parameters, the model outputs are continuous functions of the parameters. Therefore classical continuous optimization methods such as the generalized least square method can be applied (*cf.* Walter and Pronzato (1995) for a general presentation or Guo *et al.* (2006) in the context of GreenLab). Regarding topological parameters, equations (5) and (6) include transformations from real values to integers. Hence they do not generate continuous variations of the model outputs and heuristic methods need to be used. We chose the simulated annealing algorithm (Kirkpatrick *et al.* 1983). As no quantitative topological information was recorded concerning the branching structure of branches of order 2, the topological parameters were estimated using biomass data only: as a result, the optimized simulation generates average structures with global demands at each growth cycle similar to the real ones. Data of both trees were fitted simultaneously with the same set of parameter values, in order to test if the model was generic and robust enough to represent the growth of two different trees. Only the parameter *V* in Eqn (11) representing the environmental conditions was allowed to vary between the two trees.

## Results

*Direct measurements of model parameters*

Trees' ages were estimated from the observed numbers of growth units and rings. The first tree (tree 1) was 21-year-old and the second (tree 2) was 46-year-old. Several characteristics of these trees are presented in Table 1. It shows that tree 1 and tree 2 were very different for their global dimensions. For instance, their total weights were respectively 1.1 kg and 22.1 kg. Tree 1 had 45 branches of order 2 with fresh mass ranging from $1.3 \cdot 10^{-4}$ kg to 0.088 kg. Tree 2 had 81 branches of order 2 with fresh mass ranging from $6.4 \cdot 10^{-4}$ kg to 1.27 kg.

**Table 1**. *Main characteristics of the measured trees (tree 1 and tree 2). Mass refers to fresh mass.*

| Tree | 1 | 2 |
|---|---|---|
| Chronological age (years) | 21 | 46 |
| Number of order 2 branches measured | 45 | 81 |
| Height (m) | 2.80 | 8.50 |
| Total mass (kg) | 1.06 | 22.10 |
| - Trunk mass (kg) | 0.52 | 13.04 |
| - Branches mass (kg) | 0.54 | 9.06 |
| Diameter at breast height (cm) | 1.5 | 6.3 |
| Specific leaf weight (g.cm$^{-2}$) | 0.0072 | 0.0093 |
| Total leaf mass (kg) | 0.23 | 2.81 |
| Total wood mass (kg) | 0.83 | 19.29 |

Average SLW ranges from 0.0072 g.cm$^{-2}$ for tree 1 (n=20, std=0.0021) to 0.0093 g.cm$^{-2}$ (n=13, std=0.0022) for tree 2 (data not shown). No significant effect of the sample position in the tree was found. From sampling data collected, allometric ratios for biomass partitioning inside the growth units were found in coherence with the literature (Comps *et al.* 1994). The following values were taken: the average ratio between new shoots for long and short shoots is 5.25; the ratio of internode weight over leaf weight is 0.065 in a new short shoot and 0.7 in a new long shoot.

*Fitting results*

Based on the data presented in the Materials and methods section, 10 physiological parameters and 10 topological parameters (driving the numbers of metamers and axes for each PA-based category) remained to be estimated. Table 2 shows the parameter values obtained through model inversion.

**Table 2**. *Estimated parameter values after simultaneous fitting on data measured on tree 1 and tree 2. The topological parameters are estimated for each zone $Z^{ik}$. We recall that $Z^{ik}$ is the set of metamers of PA i with axillary buds of PA k; $M^{ik}$ is the number of metamers in this zone and $A^{ik}$ is the number of initiated axes in this type of zone in the whole tree. Parameters $M^{ik}$ and $A^{ik}$ are defined in equations (5) and (6). For physiological parameters, the given index refers to the equation introducing that parameter.*
*\* fixed parameters (value defined a priori)*

| Topological parameters | | | | | |
|---|---|---|---|---|---|
| Bearing PA (*i*) | Axillary PA (*k*) | Metamer number parameter $M^{ik}_1$ | $M^{ik}_2$ | Axis number parameter $A^{ik}_1$ | $A^{ik}_2$ |
| 2 | 0 | 1* | [0.4 ; 0.45] | - | - |
|  | 4 | 1* | [1.1 ; ∞ ) | 0* | [0.5 ; 0.7] |
|  | 3 | 1* | [0 ; 0.2] | 0* | [0 ; 0.2] |
|  | 2 | 0* | [0 ; 0.2] | 0* | [0 ; 0.2] |
| 3 | 0 | 1* | [1 ; 1.05] | - | - |
|  | 4 | 1* | [1.4 ; ∞ ) | 0* | [0.55 ; 0.6] |

| Physiological parameters | |
|---|---|
| Parameter $Sp_0$ in production equation [eq 11] | 0.015 |
| Exponent α in production equation [eq 11] | 0.73 |
| Sink for ring demand Pr ($g^{-1}$) [eq 7] | 2.3 |
| Exponent γ for ring demand [eq 7] | 2.95 |
| Secondary sink for ring repartition ($p_{rg}$) [eq 9] $p_{rg}$ (1)=1 (reference value for PA 1) | PA2: 0.1; PA3: 0.05; PA4: 0.01 |
| Coefficient for blade influence on ring partitioning λ [eq 9] | 0.13 |
| Environment of tree 1 [eq 11] | 0.056 |
| Environment of tree 2 [eq 11] | 0.1 |

The data of both trees were fitted in parallel with a single set of parameters. As they share the same parameter values, the observed differences between them only result from their different ages and local environments. Comparing the results found for the values of *V* for the two trees provides an indicator of the relative environmental pressures they have undergone. A smaller value of *V* was found for tree 1 (0.056 versus 0.1 for tree 2). Once the model parameters are estimated, we may retrieve by simulation the historical growth data of both trees. For example, we can evaluate the weight of tree 2 when it was 21-year-old, that is to say of the same age as tree 1. We see that tree 2 was more than twice bigger (2.54 kg versus 1.10 kg). Concerning the partitioning sub-model for rings to each metamer, the influence of the blade surface located above the considered metamer is limited. The value of the parameter driving the proportion between the two allocation modes is *λ*=0.13 (see equation 9). It means that the ring diameter profile in association with the simulated leaf surfaces does not follow the Pressler rule. In compensation, a strong effect of the branching order was found in the simulation, through the secondary sinks for ring partitioning. Indeed, the value of $p_{rg}$ is much smaller for branches (PA=2, 3, 4) than for the main stem (PA=1). It means that assimilate propagation from the leaves to the cambial sinks is nearly uniform within axes and preferentially directed towards the trunk.

The comparison between measured and fitted data is represented for total compartment weights (Fig. 4A and 4B), internode weights and diameters on the trunks (Fig. 5A and 5B) and compartment biomass on PA 2 branches (Fig. 6A, B, C, D). The results on PA 4 (not shown) are not as accurate but the biomass weights are nearly negligible compared to the other kinds of branches (compartment weight of PA 4 branches is about 20 times smaller than that of PA 3 branches). For tree 1, the simulated massed of PA 2 branches are underestimated compared to real values. The trend is reversed for PA 3 and PA 4 branches (not shown) so the total weights are globally correct (Fig. 4).

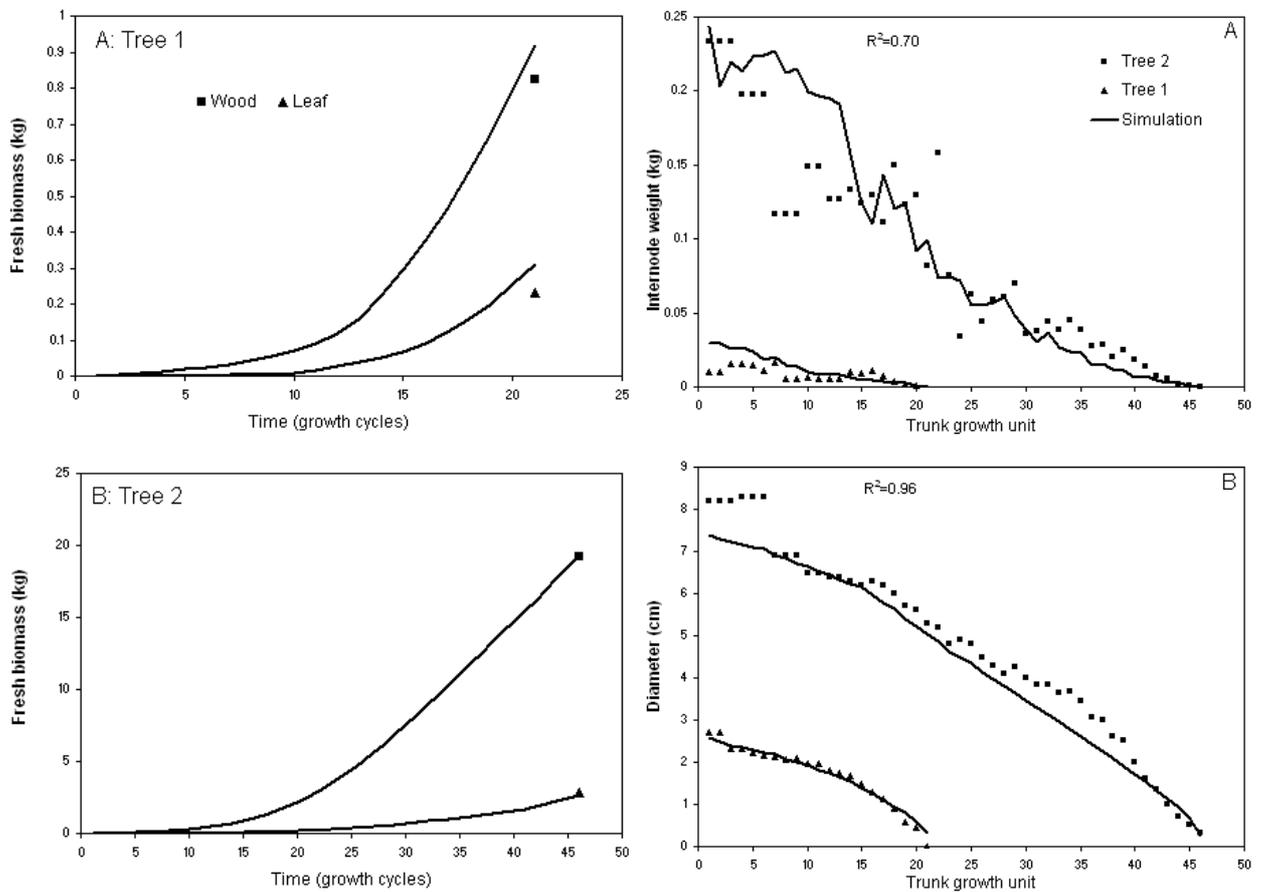

**Fig. 4.** *Evolution of biomass compartment (internodes and leaves) for measured (filled symbols) and simulated plants (lines).*

**Fig. 5**. *Average weight of one internode for each growth unit of the trunk, from base (GU 1) to top (A) and taper profile, i.e. external diameter for each growth unit of the trunk (B).*

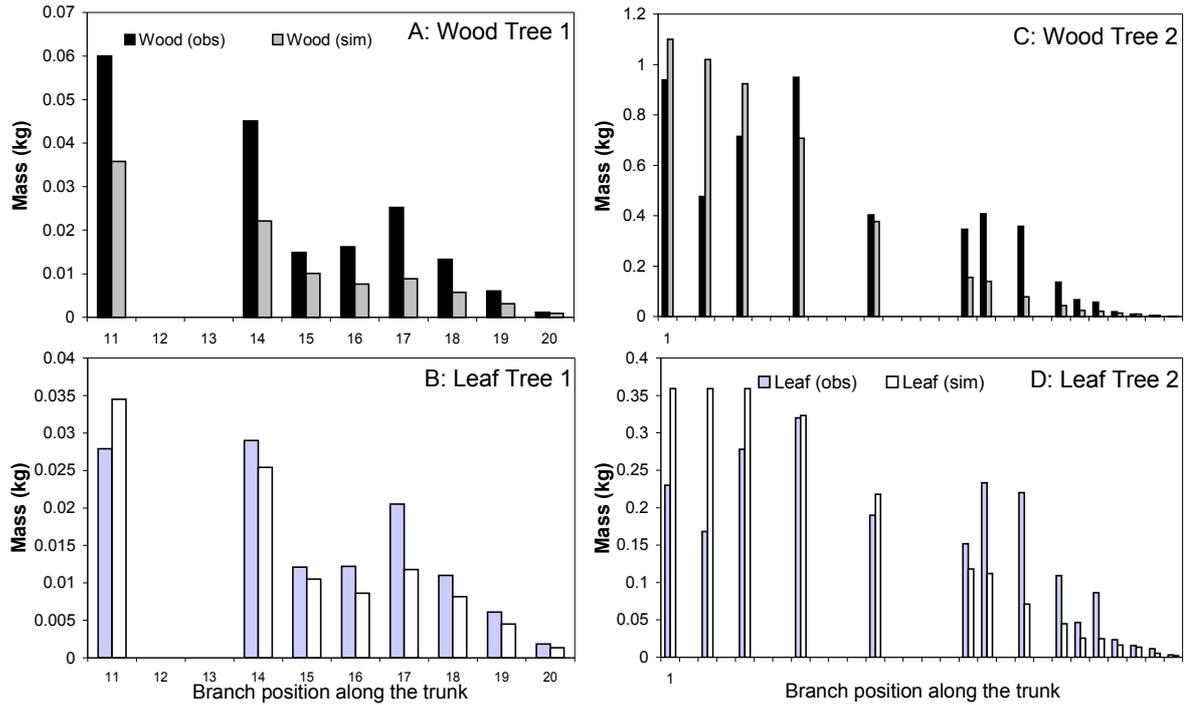

**Fig. 6**. *Fitting results: measured (black) and simulated (dark grey) wood compartment weights for branches of PA 2 (A), PA 3 (B) and PA 4 (C) for tree 1 according to the branch position on the trunk (fresh biomass). The corresponding measured (grey) and simulated (white) leaf weights are shown on the same graphs. The right graphs represent similar data for tree 2.*

Profiles of ring diameters for some growth units of tree 2 are shown in Fig. 7. The complete data set for rings is represented in a simulated-vs-measured graph in Fig. 8. The two basal growth units are not as accurately fitted as the above ones but this could be due to buttress formation which is not represented in the model. However, the global fit is satisfying ($R^2=0.92$). The accuracy of simulation of ring diameters is particularly important as it gives access to information about the past growth of the tree. Thus, not only the last stage but the complete source-sink evolution in time is likely to be correct in the simulated tree.

Note that all the graphs presented in Fig. 4 to 8 are parallel outputs of a unique fitting procedure with one single set of parameters. So the dynamics of biomass compartments presented in Fig. 4 result from the simulations of every organ mass during the plant growth.

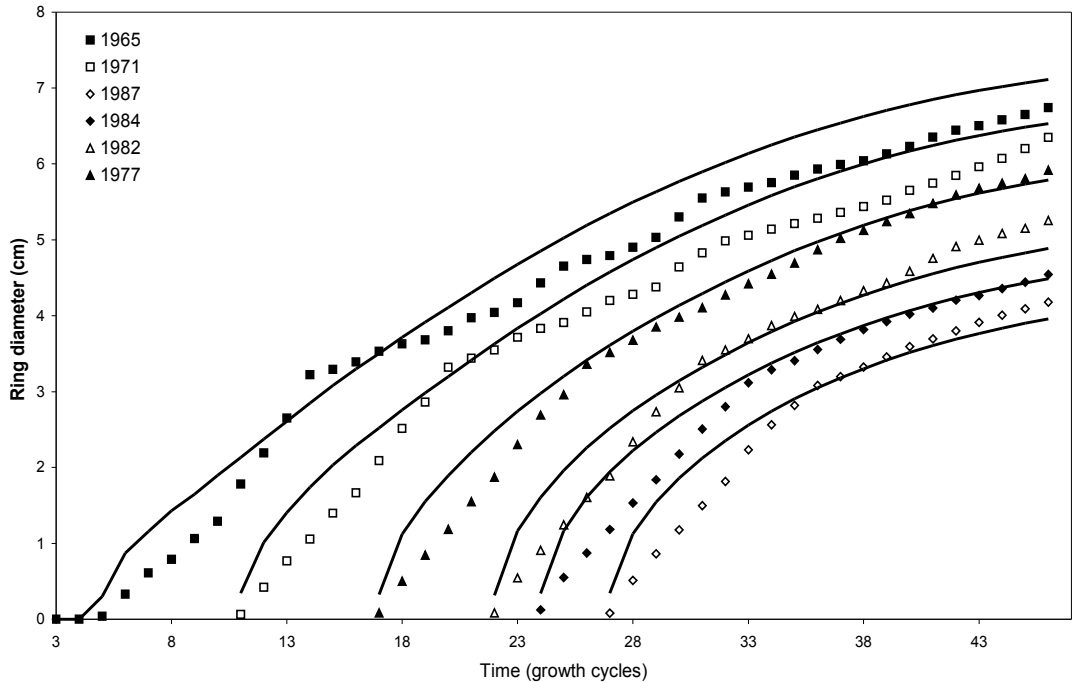

**Fig. 7**. *Diameters of some growth units of the trunk in tree 2. All the symbols of the same type represent the measured diameters of a given growth unit (named after its appearance date) as a function of the tree age. Lines are the corresponding simulated values. Each year a new ring is added. Therefore each line represents the sequence of ring diameters at a given height in the trunk.*

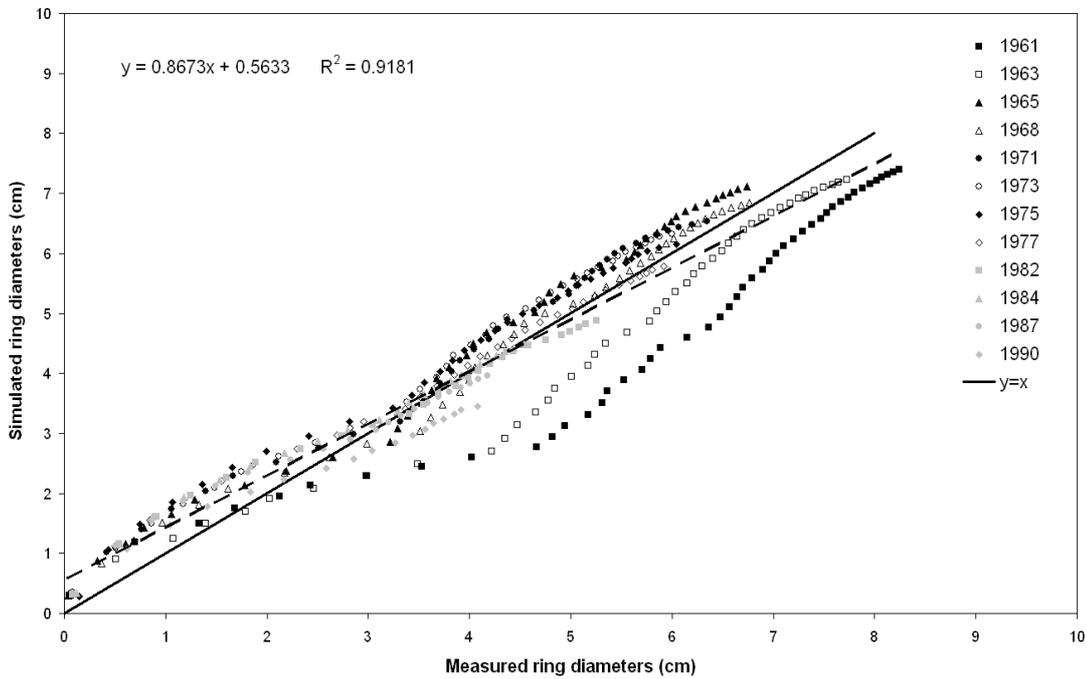

**Fig. 8.** *Observed versus simulated ring diameters for the 12 growth units (named after their appearance dates) for which the ring profiles were measured.*

Some simulation outputs from the fitted model are presented in Fig. 9 for tree 2. It shows the net biomass production at each growth cycle and its partitioning to leaves and internodes of new shoots, as well as to ring growth. A superposed curve represents the evolution of the biomass supply to demand ratio (*Q/D*, in kg). It varies in the interval [$1 \times 10^{-3}$; $4.31 \times 10^{-3}$] for tree 2 and in the interval [$9.6 \times 10^{-2}$; $2.4 \times 10^{-3}$] for tree 1. Note that *Q/D* reaches a maximum then declines and so does the biomass amount allocated to cambial growth. Concerning the topological parameters, an interval is given rather than a single value. Indeed, as equations (5) and (6) include transformations from real values to integers, the simulated structures are the same when the topological parameters take their values in these intervals. For some parameters, the upper boundary can be equal to infinity when they take their maximal values ($M^{ik}_{max}$). Getting more accuracy would require collecting data on trees grown in even more contrasted conditions (so that a larger range of variations of the *Q/D* ratio could help discriminate between the threshold values more precisely) and incorporating more botanical knowledge about tree architecture in the model, as discussed in the Discussion section.

The topological structures of the fitted trees are shown in Fig. 10 and the 3D output for tree 2 is drawn. As no topological data were recorded on branches, the number of axes cannot be precisely compared but their demands are similar enough to reproduce the biomass allocation to each compartment and to each growth unit of the stem. Note that for each branch of order 2, there is a transitory phase of installation of the branching structure: there are progressive increases of number of metamers, number of axes and vigor of axes along the branch development (see for instance PA 2 axis of tree 1 in Fig. 10). That basal effect is not forced but dynamically generated by the model thanks to the feedback effect driven by the *Q/D* variable.

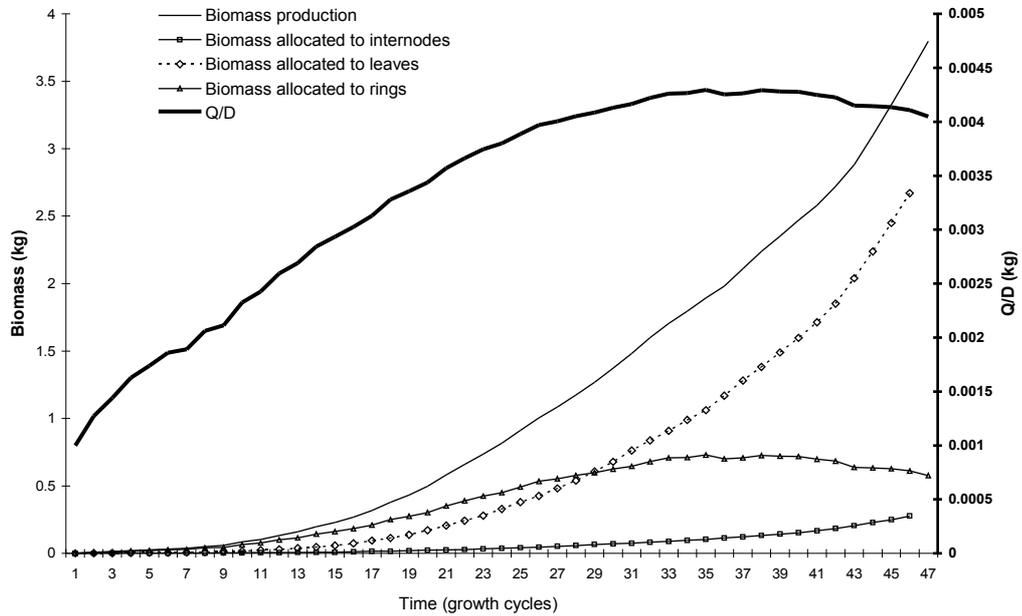

**Fig. 9**. *Simulation outputs for tree 2: net biomass production and partitioning at each growth cycle (left Y-axis); evolution of the ratio of biomass production to demand, Q/D (right Y-axis).*

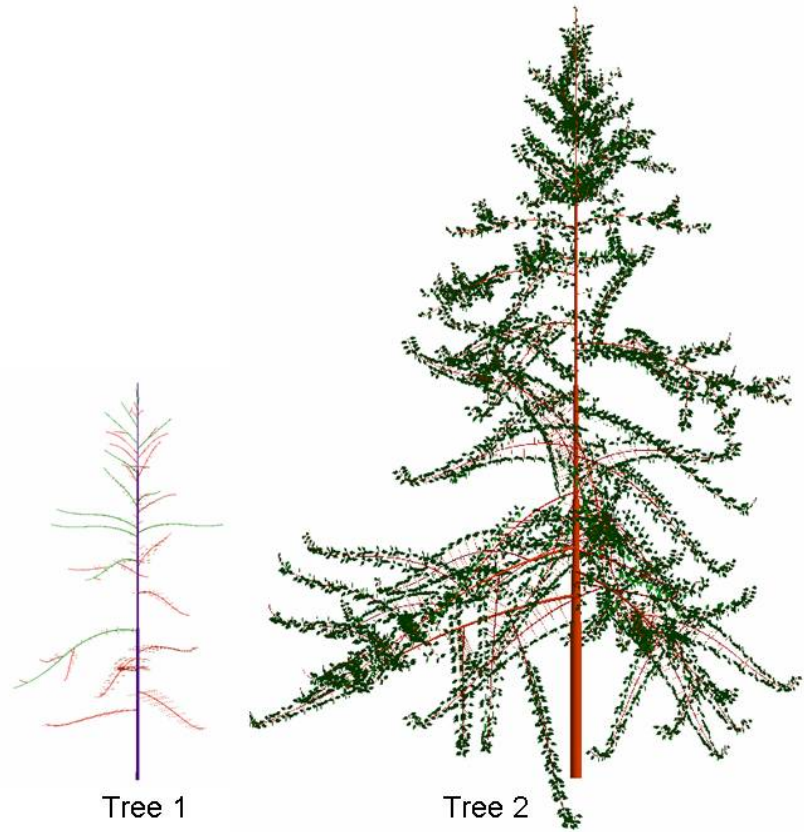

**Fig. 10**. *Topological structure of tree 1 and 3D output of tree 2 resulting from simulations with the estimated model parameters. The number of branches and metamers per growth unit are driven by the ratio of biomass supply to demand (Q/D).*

## Discussion

This preliminary study shows that a functional-structural growth model based on a simplified description of tree structure can mimic the morphological plasticity and the source-sink dynamics of two understorey beech trees at different growth stages with a single set of species-specific parameters. The architectural structure is considered controlled by a dynamic feedback of tree internal trophic state on the organogenesis processes and factorized in order to take advantage of the repetitions in tree architecture. It results in a natural simplification of the topology: it is no longer based on a detailed description of each branch growth and development, which would entail facing an inextricable variability, but the global source-sink dynamics thus generated are close enough to reproduce the biomass allocation to each compartment and to each growth unit of the stem. The model parameters were estimated by model inversion with a set of experimental data in adequacy with the simulation outputs. This work is a preliminary study on the feasibility of such procedure, and further investigation is needed in order to give a predictive value to the model.

Besides the simplified representation of topology, several other simplifications require discussion. Firstly, allocation to the root system was considered as a constant proportion of the total biomass produced. This could be modified by defining the sink evolution of the root systems if the corresponding data were available.

The variations of specific leaf weight according to leaf position in the crown were not considered, although several studies have revealed that a gradient can be observed in relation to light exposure (Beaudet and Messier 1998; Sarijeva *et al.* 2007) and to competitive status (Barna 2004). However, no visible gradient was observed from our sampling data. As our trees have grown in the understorey and were shaded by higher neighbor trees, the light they received was mainly diffuse light. So the light environment was likely to be relatively isotropic and similar for both trees

Since no precise data were registered concerning the inter-annual climate variations, it was chosen to input a constant value to represent an average environmental variable, distinct for each tree. As a consequence, simulations mainly represent the ontogenetic trends induced by the dynamics of trophic competition. If environmental data were available, the effects of inter-annual climatic variations would also be reflected in the source-sink dynamics described by the model. A lower value was found for the environmental potential of tree 1 than of tree 2. This is consistent with the qualitative observations on the field: tree 1 was surrounded by close neighbor trees and therefore had probably suffered from higher competition than tree 2 that was grown in a more open space.

The value of the exponent $\gamma$ for ring demand is found to be nearly 3. It is coherent with the common observation that ring width is particularly sensitive to changes of growth conditions. Since the environmental conditions were considered as constant in simulations, $\gamma$ is likely to be over-estimated to account for higher variations than those generated by source-sink dynamics (namely the ratio *Q/D*).

Beech trees have propensity to develop epicormic branches (emerged from latent or adventitious buds, with possibly high delays after the bud formation) that lower timber quality. It occurs in particular growth phases and at specific locations in trees, in response to abrupt changes in tree growth conditions or to stresses as an attempt to replace the deficient crown (Nicolini *et al.* 2001). That phenomenon was observed on tree 2 but was not modeled, as the biomass weights involved were negligible. However, it could be interesting to study its interaction with the model variables. Nicolini *et al.* (2001) showed that appearance of epicormic branches was associated to low cambial activity for understorey beech trees. This is coherent with our results that reveal a stabilization and a decrease in the Q/D ratio from the 32$^{nd}$ growth cycle that coincides approximately with the emergence of short epicormic axes on the basal growth units of the trunk (around the 30$^{th}$ growth cycles).

If these first results on the parametric identification of GreenLab on adult trees are encouraging, further research work are necessary to improve this procedure, in particular by incorporating more botanical information to refine the values of the topological parameters. Some potential improvements are discussed below.

First, a refinement of the fitting methodology introduced in this paper would be to include topological criteria as target data. Based on the formalism introduced by Godin and Caraglio (1998) that represents the tree as a rooted multi-scale graph, methods have been developed to compare tree architectures (Ferraro and Godin 2000). Thus the distance between two tree structures can be incorporated in the definition of the objective function to minimize. It could provide a precious indicator about the reliability of the simulated architecture. Ideally, the distance between the whole crown architectures of the simulated and the real trees would be taken into account. However, in our modeling approach, the simplifications imposed by the complex structure of a tree crown imply that the exact replication of its topology is nearly impossible. So the distance should be computed at an intermediate level to compare structures of branch samples from the target tree to the equivalent average structure of the simulated tree.

Moreover, crown architecture is the result of a quasi-infinite number of interactions between unknown events (weather, local competition, insect attacks…) generating the stochastic aspect of structure, especially for beech trees. Simulations cannot retrace the historic succession of events and is thus intrinsically restricted to reproducing only its average or global dynamics. But botanical studies from exhaustive or sampling measurements of tree architecture have revealed species-specific and ontogenetic patterns. It means that in branching structures, some temporal and/or spatial trends can be identified; they are characteristic of species and largely independent from the environmental conditions. These particular trends can be quantified according to adequate variables such as distributions of shoot lengths. It can be represented through the parameters of hidden semi Markov chains (Guédon *et al.* 2001) or hidden Markov tree models (Durand *et al.* 2005). That botanical knowledge can be used to build realistic tree mock-ups. It provides an elegant and compact way of describing and characterizing a branching architecture. But it still needs to be linked with a functional part if issues related to yield prediction are needed. Our long-term objective is to relate the parameters of those botanical descriptions to functional processes and environmental controls. That kind of work was initiated in the L-Peach model (Allen *et al.* 2005): bud fates in each zone of a shoot are regulated by the quantity of carbon available, using hidden semi-Markov chains to reproduce the observed stochastic patterns of Peach tree branching systems (Costes *et al*. 2007). The results are qualitatively satisfying but the authors underline that more experimental work is needed for the model calibration, especially to investigate the values of carbon availability threshold for bud breakouts in each zone. In our opinion, the work presented in this paper is complementary to their approach as it exposes a possible way to extract some threshold values from experimental data.

Finally, the work reported in this paper can be seen a preliminary step towards the integration of the knowledge and techniques acquired in the field of functional-structural modeling to simulate the growth of individual-based spatialized stands on long periods and taking into account architectural information, for forestry applications.

**Acknowledgement**
Data were collected with the support of the LERFOB Wood Quality research team. We particularly thank Claude Houssement, Emmanuel Cornu and Alain Mercanti for their help in the measurement process. This research is partly supported by the Cap Digital Business Cluster **Terra Numerica** project.